\documentclass[11pt,a4paper]{article}

\usepackage[english]{babel}
\usepackage{latexsym}
\usepackage{graphicx}
\usepackage{enumerate}

\usepackage{amsthm,amsmath,amsfonts,amssymb}

\newcommand\blfootnote[1]{%
	\begingroup
	\renewcommand\thefootnote{}\footnote{#1}%
	\addtocounter{footnote}{-1}%
	\endgroup
}

\def\A{\mbox{\boldmath $A$}}

\def\I{\mbox{\boldmath $I$}}
\def\J{\mbox{\boldmath $J$}}

\def\dim{\mathop{\rm dim}\nolimits}
\def\dist{\mathop{\rm dist}\nolimits}

\def\ker{\mathop{\rm ker}\nolimits}
\def\mod{\mathop{\rm mod}\nolimits}
\def\>{\mathop{\rightarrow}\nolimits}

\def\i{\mbox{\boldmath $i$}}
\def\j{\mbox{\boldmath $j$}}
\def\m{\mbox{\boldmath $m$}}

\def\v{\mbox{\boldmath $v$}}

\def\vec0{\mbox{\bf 0}}

\def\sp{\mathop{\rm sp}\nolimits}

\def\0{\mbox{\bf 0}}
\def\1{\mbox{\bf 1}}

\begin{document}

\title{Graphs, friends and acquaintances}

\author{C. Dalf\'o$^a$, M.A. Fiol$^b$\\
$^{a}${\small Departament de Matem\`atica} \\
{\small Universitat de Lleida} \\
{\small Igualada (Barcelona), Catalonia} \\
{\small {\tt{cristina.dalfo@matematica.udl.cat}}} \\
$^{b}${\small Departament de Matem\`atiques} \\
{\small Universitat Polit\`ecnica de Catalunya} \\
{\small Barcelona Graduate School of Mathematics} \\
{\small Barcelona, Catalonia} \\
{\small{\tt{miguel.angel.fiol@upc.edu}}}}

\date{}

\maketitle


\vskip0.3cm


\begin{abstract}
As is well known, a graph is a mathematical
object modeling the existence of a certain relation between pairs of
elements of a given set. Therefore, it is not surprising that many of the first results concerning graphs made reference to
relationships between people or groups of people. In this article,
we comment on four results of this kind, which are related to various
general theories on graphs and their applications:
the Handshake lemma (related to graph colorings and Boolean algebra), a lemma on known
and unknown people at a cocktail party (to Ramsey theory), a
theorem on friends in common (to distance-regularity and coding theory), and Hall's Marriage theorem (to the theory of networks). These four areas of graph theory, often with problems which are easy to state but difficult to solve, are extensively developed and currently give rise to much research work. As examples of representative problems and results of these areas, which are discussed in this paper, we may cite the following: the Four Colors Theorem (4CTC), the Ramsey numbers, problems of the existence of distance-regular graphs and completely regular codes, and finally the study of topological proprieties of interconnection networks.
\end{abstract}

\blfootnote{
	\begin{minipage}[l]{0.3\textwidth} \includegraphics[trim=10cm 6cm 10cm 5cm,clip,scale=0.15]{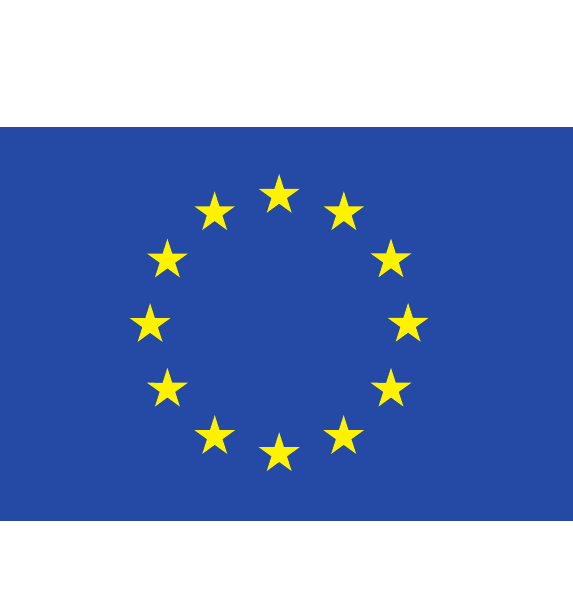} \end{minipage}  \hspace{-2cm} \begin{minipage}[l][1cm]{0.79\textwidth}
		The research of C. Dalf\'{o} has also received funding from the European Union's Horizon 2020 research and innovation programme under the Marie Sk\l{}odowska-Curie grant agreement No 734922.
\end{minipage}}

\section{Introduction}

A graph $G=(V,E)$ is a mathematical structure consisting of a
vertex set $V$ and a set of edges $E$ (or nonordered pairs of
vertices). Normally, each vertex $v\in V$ is represented by a point
and each edge $e=\{u,v\}\in E$ by a line joining vertices $u$ and $v$. Graph
theory belongs to combinatorics, which is the part of
mathematics that studies the structure and enumeration of discrete objects, in contrast to the continuous objects studied in mathematical analysis.
In particular, graph theory is useful for studying any system with a certain relationship between pairs of elements, which give a binary relation. It is therefore
not surprising that many of the problems and results were originally stated
in terms of personal relationships. For example, one of the most simple
results is the Handshake lemma: \emph{At a cocktail party,
an even number of people shake an odd number of hands}.
There is also the so-called Friendship
theorem: \emph{At a party, if each pair of people has
exactly one friend in common, then there is somebody who is friend of
everybody}. The first and most appealing proof of this theorem is
due to Paul Erd\H{o}s (with Alfred R\'{e}nyi and Vera S\'{o}s), a
Hungarian mathematician, probably the most prolific of the 20th
century, who like Euler enjoyed coining sentences such as ``A mathematician
is a device for turning coffee into theorems'' or ``Another roof,
another proof''. The latter phrase shows his great capacity and
predisposition for collaborating with other authors from all over the world (he had 509 coauthors). From Erd\H{o}s we have the \emph{Erd\H{o}s number}: the co-authors of Erd\H{o}s have Erd\H{o}s number 1, the co-authors of the co-authors of Erd\H{o}s have Erd\H{o}s number 2, etc. For more information on Erd\H{o}s, see Hoffman~\cite{Ho98}.

It is considered that the first paper on graph theory was published in 1736. Its
author was the great Swiss mathematician Leonhard Euler, about who
it is said that he wrote papers in the half an hour between the first
and the second calls for lunch. This first paper is about the existence of a possible walk
across the K\"{o}nigsberg bridges; see Euler~\cite{e1741}. This
city was the capital of Oriental Prussia, the birthplace of Immanuel Kant.
Nowadays it corresponds to the Russian city of Kaliningrad.
The problem of the K\"{o}nigsberg bridges is related to the puzzle
of drawing a figure without raising the pencil from the paper and
without passing twice through the same place. In the original
problem, it was asked if it was possible to walk through the city by
crossing all the bridges only once. With an ingenious reasoning, which in fact does not explicitly use any graph, Euler proved the impossibility of this walk.

Another of the most famous problems in graph theory, not solved
until 1977 by Appel, Haken and Kock~\cite{ApHaKo77,Ap77}, is the Four Colors theorem (4CT), which states that the countries of any map drawn in the plane can be colored with four
colors, such that countries with a common border (different from a
point) bear different colors. This theorem is regarded as the first important result to be proved using a computer, because in a part of its proof
1,482 configurations were analyzed. For this reason, not all mathematicians accept it. Twenty years later,
Robertson, Sanders, Seymour and Thomas~\cite{RoSaSeTh97} gave an
independent proof, which is shorter, but also requires the use of
a computer, because of the 633 configurations analyzed.

As we have already stated, graph theory is used to study different relations.
A first example is an electric circuit, with all its components and its
connections. In telecommunications, graph theory
contributes to the modeling, design and study of interconnection or
communication networks. For instance, interconnection networks are
used in multiprocessor systems, where some processors undertake
a task of exchanging information, and in local networks consisting of
different computers placed at a short distances, which exchange data at very
high speed and low cost. As regards communication networks, nowadays the most
important example is the Internet, which makes the
communication and exchange of data possible between computers all
around the world. In fact, we are experiencing a communication
revolution, so that we could say that we are `weaving' the
communication network.

For more details about notation, basic concepts and history of graph
theory see, for example, Bollob\'{a}s~\cite{Bo90}, Diestel~\cite{Di97}, West~\cite{w00} and Biggs, Lloyd and
Wilson~\cite{blw76}.

\section{Shaking hands: Colorings and Boolean algebra}

In a graph $G=(V,E)$, the \emph{degree} $\delta(u)$ is the number of adjacent vertices to vertex $u$,
namely, the number of incident edges to $u$. We denote by $\Delta(G)$ the maximum degree of all the vertices of $G$ and by $\delta(G)$ the minimum degree.

We begin with one of the most simple results about graphs, which
states that the sum of the degrees of the vertices in $V$ equals
twice the number of edges in $E$:
\begin{equation}
\label{encaixada} \sum_{u \in V}\delta (u) = 2 |E|,
\end{equation}
since in the degree sum, we count each edge twice
because each edge is incident to two vertices.
From here, we obtain the inequalities:
\begin{equation}
\label{desigualtat} \delta(G)|V|\leq2|E|\leq\Delta(G)|V|.
\end{equation}
Although these
results are apparently trivial, they have some interesting
corollaries, such as the following:
\begin{itemize}
\item[$(a)$]
{\it Every graph has an even number of vertices with odd degree}.

This is the so-called Handshake lemma, because it can be stated
as follows: {\it At a cocktail party,
the number of people who shake an odd number of other people's hands is always even}.
\item[$(b)$]
{\it Every $\delta$-regular graph $($a graph is $\delta$\emph{-regular} if all its vertices have
degree~$\delta)$, with $\delta$ odd, has an even number of
vertices}.
\item[$(c)$]
{\it Every planar graph $($that is, it can be drawn on the plane without edge crossings$)$
with girth $g$ $($the \emph{girth} is the length of the
shortest cycle$)$ and number of edges $|E|$ satisfies}
\begin{equation}
\label{eq-nombrearestes}
 |E| \le \frac {g (|V|-2)}{g-2}.
\end{equation}
\end{itemize}

To prove $(c)$, we need the well-known Euler formula~\cite{e1752}
published between 1752 and 1753, and already observed by Descartes in 1640,
which can be proved by induction and states that every planar graph with $n=|V|$ vertices, $m=|E|$
edges and $r=|R|$ regions satisfies
\begin{equation}\label{eq-Euler}
r+n=m+2.
\end{equation}

\begin{figure}[t]
\begin{center}
\includegraphics[width=12cm]{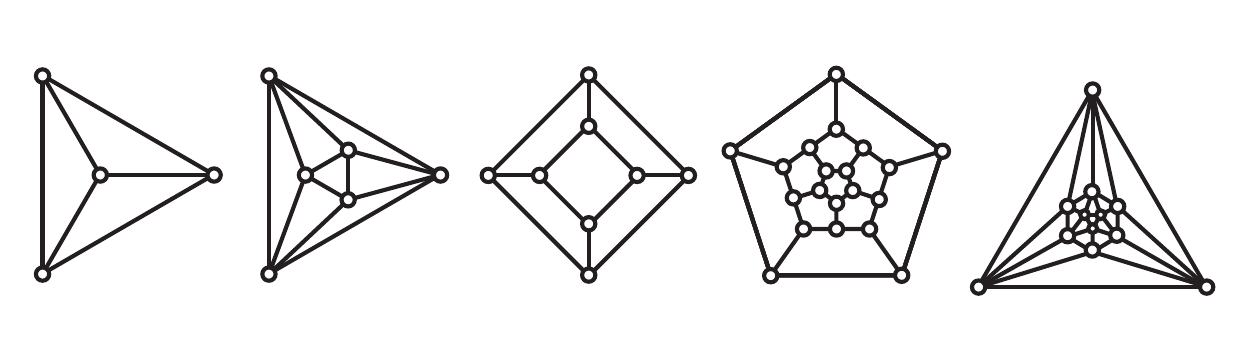}
\caption{The graphs of the five Platonic solids.}
\label{fig:solids_platonics}
\end{center}
\end{figure}

In this formula, the number of regions includes
the exterior one (that is, the `sea', if we have a map or if the
graph is imbedded on a sphere). For example, the Euler formula is
satisfied by the graphs of the Platonic solids shown in
Figure~\ref{fig:solids_platonics}. In fact, this formula gives
necessary conditions for the existence of these regular polyhedra;
see Rademacher and Toeplitz~\cite{RaTo70}.
In proving (\ref{eq-Euler}), the key fact is that the removing of a vertex with degree $\delta$ (and its incident edges) leaves a new planar graph whose number of regions, vertices and edges have been reduced, respectively, by $\delta-1$, $1$ and $\delta$ units.

\begin{figure}[t]
\begin{center}
\includegraphics[width=5cm]{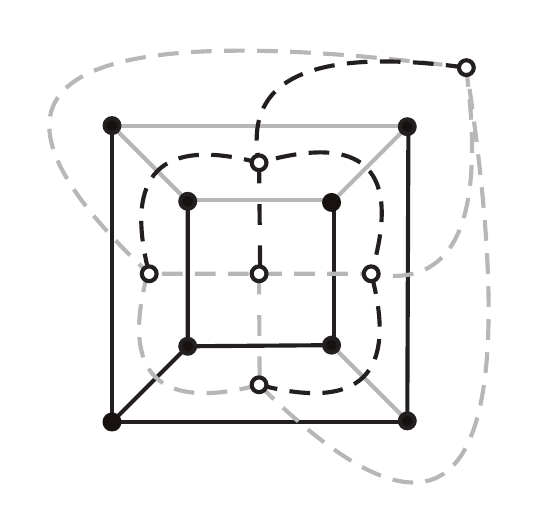}
\caption{The spanning tree (black edges) of the cube graph $Q$
(continuous edges and black vertices) and its dual (dashed edges and
white vertices).} \label{fig:2arbres_generadors}
\end{center}
\end{figure}

Returning again to the Euler formula, the number $r$ of regions
can also be interpreted as the cardinality of the vertex set of the dual graph $G^*$. Given a planar
graph $G$ with $n=|V|$ vertices and $m=|E|$ edges forming regions,
its \emph{dual} graph $G^*=(V^*,E^*)$ has vertices representing the regions of $G$,
and there is an edge between two vertices if the corresponding
regions are neighbors. Then, $r=|V^*|$ and $m=|E|=|E^*|$. This interpretation provides a more
symmetric Euler formula:
\begin{equation}
\label{eq.euler}
|E^*| = (|V^*|-1)+(|V|-1) = |E|,
\end{equation}
which allows us to prove it without using induction, but rather by
identifying both parenthesis in Equation (\ref{eq.euler}) as the
number of edges of two spanning trees $T^*$ and $T$ belonging to
$G^*$ and $G$,
respectively. A \emph{spanning tree} $T$ of a \emph{connected} graph
$G=(V,E)$ (that is, there is a path between any pair of vertices)
is composed of the vertex set $V$ and $|V|-1$ edges without forming cycles.
An example of this is shown in Figure~\ref{fig:2arbres_generadors}, where each black
continuous edge of $G$ (the graph of a cube $Q$) belongs to $T$, but where
each black dashed edge corresponds to an edge
of $T^*$ in $G^*$ (the graph of an octahedron). For more details, see
Aigner and Ziegler~\cite{az98}.

In our case, the proof of $(c)$ is as follows: As each
edge is the border of two regions and each region has at least $g$
edges, we have $r\leq 2m/g$.
Note that this inequality is obtained from (\ref{desigualtat}), considering the dual graph, since $r=|V^*|$, $m=|E^*|$ and $g=\delta(G^*)$.
Using this inequality and Equation~$(\ref{eq-Euler})$,
we obtain $(\ref{eq-nombrearestes})$.

\begin{figure}[t]
\begin{center}
\includegraphics[width=7cm]{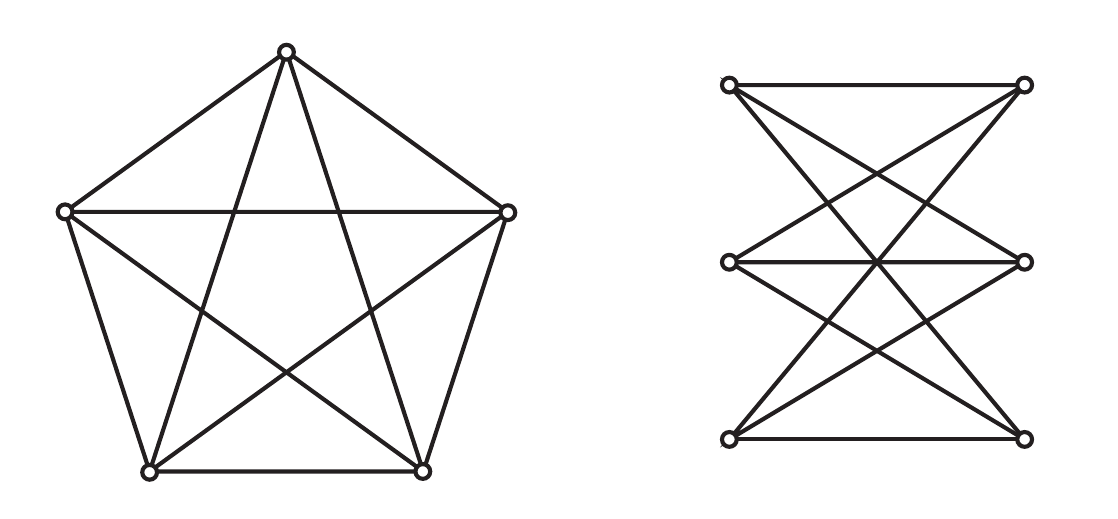}
\caption{The complete graph $K_5$ and the complete bipartite graph
$K_{3,3}$.} \label{fig:K5-K33}
\end{center}
\end{figure}

As a particular case of $(c)$, we have the following
result:
\begin{itemize}
  \item[$(d)$] \emph{In any planar graph $(g\geq3)$ the number of edges satisfies $m \le 3n-6$;
  if it does not contain triangles $(g\geq4)$, then $m\le 2n-4$; and
  if it  contains neither triangles nor squares $(g\geq5)$, then $m\le \frac{5}{3}(n-2)$.}
\end{itemize}
From the first inequality, we can see that the complete graph $K_5$
$(n=5$, $m=10)$ is not planar. A graph is \emph{complete} if there is an edge between every pair of vertices.
Similarly, from the second inequality, we also obtain that
the complete bipartite graph $K_{3,3}$ $(n=6,m=9)$ is not planar.
A \emph{bipartite} graph (that is, the vertex set can be decomposed into two independent subsets such that vertices in every subset are not adjacent) is \emph{complete} if each pair of vertices in different subsets are adjacent.
See both graphs in Figure~\ref{fig:K5-K33}. Notice that, for instance, the third inequality
turns out to be an equality in the case of the dodecahedron graph
(see again Figure~\ref{fig:solids_platonics}, $n=20$ and
$m=30$).

In this context, we have the famous Kuratowski theorem~\cite{Ku30}, which characterizes planar graphs (see also the book by West~\cite[pp. 246--251]{w00} and the paper by Thomassen~\cite{Tho90}, where the relation between the planarity criterion and the Jordan Curve Theorem is explained):
\begin{itemize}
\item \emph{A graph is planar if and only if it contains no homeomorphic subgraph to $K_5$ or $K_{3,3}$.}
\end{itemize}
Recall that a graph $H$ is {\em homeomorphic}
to a graph $G$ if the edges of $G$ correspond to (independent) paths in $H$.

From Equation (\ref{encaixada}) and again the inequalities in $(d)$, we can prove the following:
\begin{itemize}
  \item
  \emph{Every planar graph $G$ contains a vertex $u$ of degree
$\delta(u)\le5$. Moreover, if $G$ does not contain triangles, then it has a vertex $u$ of degree \linebreak $\delta(u)\le3$.}
\end{itemize}

Indeed, if $n_i$ denotes the number of vertices with
degree $i\in\mathbb{N}$, then from Equation (\ref{encaixada}) we have that
$$
2m=n_1+2n_2+3n_3+\cdots\leq2(3n-6)=6n_1+6n_2+6n_3+\cdots-12,
$$
whence
$$
5n_1+4n_2+3n_3+2n_4+n_5-n_7-2n_8-\cdots=12,
$$
so that $n_i\geq0$ for some $i\leq5$, as claimed. The proof of the case without triangles is analogue.

The existence of a vertex with degree at most five allows us to prove, by induction, the Five Color
theorem (5CT), which was first proved by Heawood~\cite{h1890}
(see Aigner and Ziegler~\cite{az98}):
\begin{itemize}
  \item \emph{Five colors suffice to get a vertex-coloring of a planar graph.}
\end{itemize}
Recall that in a vertex-coloring, adjacent vertices have different
colors.

\begin{figure}[t]
\begin{center}
\includegraphics[width=6cm]{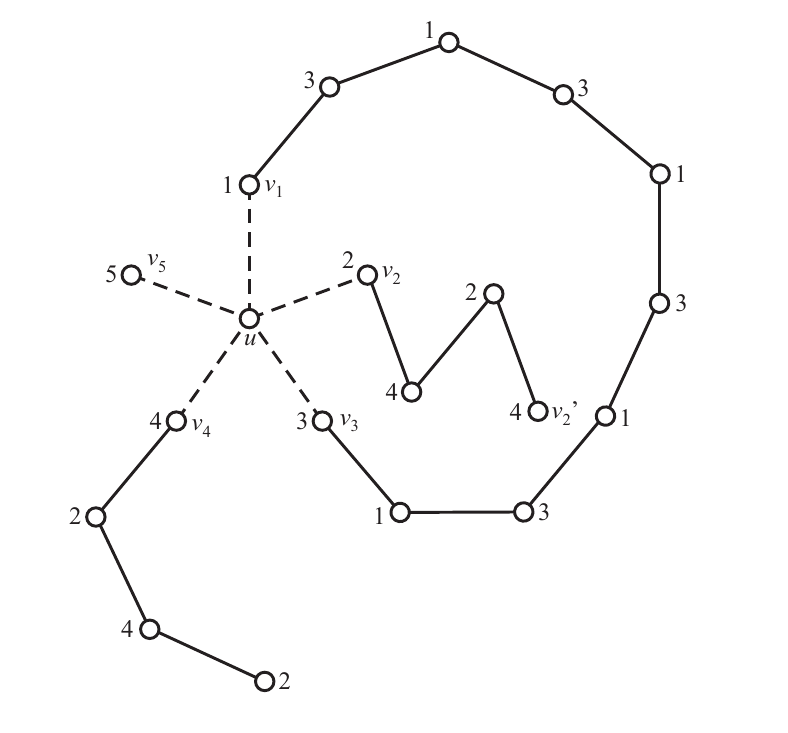}
\caption{The case $r=\delta=5$ in the proof of the Five Color theorem (5CT).}
\label{fig:Heawood}
\end{center}
\end{figure}

First note that the result is trivially true for graphs with at most 5 vertices. Then, assume that it is also true for graphs with $n-1>5$ vertices, and let $G$ be a graph with $n$ vertices. We know that $G$ contains a vertex $u\in V$ with degree $\delta\le 5$. Let $v_i$, $1\le i\le\delta$, denote the adjacent vertices to $u$. From the induction hypothesis, the graph $G'=G-u$ (obtained from $G$ by removing vertex $u$ and all its incident edges) has a vertex-coloring with $r\le 5$ colors.
Therefore, if $r\le 4$ (which is always the case when $\delta\le 4$), we can restore vertex $u$ and give it a color different from the colors of the adjacent vertices $v_i$. Thus, we obtain a coloring of $G$ using at most 5 colors.
Otherwise, if $r=\delta=5$ we can assume, without lost of generality, that we have a situation as shown in Figure~\ref{fig:Heawood} (where vertex $v_i$ has color $i$, $1\le i\le \delta$). Now consider the paths with vertices alternatively colored 1-3 (with final vertices $v_1$ and/or $v_3$) and 2-4 (with final vertices $v_2$ and/or $v_4$). As $G'$ is planar, these possible paths cannot cross each other (that is, they have neither crossed edges nor common vertices). Then if, for example, there exists the path 1-3 with initial-final vertices $v_1$-$v_3$, the path 2-4 with initial vertex $v_2$ cannot have $v_4$ as final vertex, but another vertex denoted by $v_2'$ (see again Figure~\ref{fig:Heawood}). Therefore, we can interchange the colors 2-4 in this path, so that $v_2$ gets color 4. We can then restore vertex $u$ and assign it color $2$, obtaining a coloring of $G$ with 5 colors.

We now consider the case of giving one of three
colors to each edge of a graph $G$ with maximum degree $3$. This
is called a \emph{free edge-coloring} of $G$. In
particular, the (`not-free') \emph{edge-coloring} of a cubic ($3$-regular)
graph, also called {\it Tait-coloring}, corresponds to the case
where adjacent edges receive different colors. As we will see
later, if $G$ is a planar graph, the problem of the existence of
Tait-colorings is closely related to the Four Color theorem (4CT).
Moreover, we will also see that the construction of cubic graphs which
cannot be Tait-colored leads to Boolean algebra, which is commonly
used in the study of logic circuits. To this end, we introduce a natural
generalization of the concept of `color', which describes in a
simple way the coloring (``$\0$'' or ``$\1$'') of any set of edges
or, more abstractly, of any family ${\cal F}$ of $m$ colors chosen between
three different colors, say ${\cal C}=\{1,2,3\}$, such that color $i\in {\cal C}$ appears
$m_i$ times. This situation can be represented by the coloring-vector
$\m=(m_1,m_2,m_3)$, where $m=m_1+m_2+m_3$. Then, we say that $\cal F$ has {\it
Boole-coloring} $\0$, denoted by $\Psi({\cal F})=\0$, if
$$
m_1\equiv m_2\equiv m_3 \equiv m \quad (\mod 2),
$$
whereas $\cal F$ has {\it Boole-coloring} $\1$ (more specifically $\1_a$), denoted by  $\Psi({\cal F})=\1$ (or $\Psi({\cal F})=\1_a$),
if
$$
m_a+1\equiv m_b\equiv m_c\equiv m+1 \quad (\mod 2),
$$
where $\{a,b,c\}=\{1,2,3\}$. See Fiol and Fiol~\cite{f&f84} for more information.

Recalling these definitions, the Boole-coloring of an edge $e\in E$
with color $a\in\cal{C}$ is $\Psi(e)=\Psi(\{a\})=\1_a$, and the
Boole-coloring of a vertex $v\in V$, denoted by $\Psi(v)$, is defined
as the Boole-coloring of its incident edges, which can have either different
or the same colors. In this context, it is curious to note the
following facts:
\begin{enumerate}
\item
If $\delta (v)=1$, then $\Psi (v)=\1_{a}$ if and only if the incident edge
to vertex $v$ has color $a\in\cal{C}$.
\item
If $\delta (v)=2$, then $\Psi (v)=\0$ if both incident edges to
vertex $v$ have the same color, and $\Psi (v)=\1$ if not.
\item
If $\delta(v)=3$, then $\Psi (v)=\0$ if and only if the three incident edges
to vertex $v$ have three different colors. Thus, in a
Tait-coloring of a cubic graph, all its vertices have Boole-coloring $\vec0$.
\end{enumerate}

Moreover, a natural sum operation can be defined in the set ${\cal
B}=\{\0,\1_{1}$, $\1_{2},\1_{3}\}$ of Boole-colorings in the following
way: Given the colorings $X_1$ and $X_2$ represented, respectively, by the
coloring-vectors $\m_1=(m_{11},m_{12},m_{13})$ and
$\m_2=(m_{21},m_{22},m_{23})$, we define the sum $X=X_1+X_2$ as the
coloring represented by the coloring vector $\m=\m_1+\m_2$. Then,
$(\cal{B},+)$ is isomorphic to the Klein group, with $\0$ as
identity, $\1_{a}+\1_{a}=\0$, and $\1_{a}+\1_{b}=\1_{c}$ where $\{a,b,c\}=\{1,2,3\}$; see
Table~.

\begin{table}
\begin{center}
\begin{tabular}{|c||c|c|c|c|}
\hline
$+$ & $\vec0$ & $\1_1$ & $\1_2$ & $\1_3$\\
\hline\hline
$\vec0$ & $\vec0$ & $\1_1$ & $\1_2$ & $\1_3$\\
\hline
$\1_1$ & $\1_1$ & $\vec0$ & $\1_3$ & $\1_2$\\
\hline
$\1_2$ & $\1_2$ & $\1_3$ & $\vec0$ & $\1_1$\\
\hline
$\1_3$ & $\1_3$ & $\1_2$ & $\1_1$ & $\vec0$\\
\hline
\end{tabular}
\caption{Klein's group of Boole-colorings.}
\end{center}
\label{taula:klein}
\end{table}

Notice that, since every element coincides with its inverse,
$m\1_{a}=\1_{a}+\1_{a}+\stackrel{m}{\cdots}+\1_{a}$ is $\0$ if $m$
is even and $\1_{a}$ if $m$ is odd. From this simple fact, we can
imply the following result (see Fiol~\cite{f95}), which is very useful in the further
development of the theory and can be regarded as a generalization
of the so-called Parity lemma (see Isaacs~\cite{i75}):

\begin{itemize}
\item {\it Let $G$ be a graph with $n$ vertices, maximum degree $3$, and having a
free edge-coloring, such that $n_{i}$
vertices have Boole-coloring $\1_{i}$, for $i\in \cal{C}$, with $n'=n_1+n_2+n_3\le n$. Then,
\begin{equation}\label{parity-lemma}
n_{1} \equiv  n_{2} \equiv  n_{3} \equiv  n'\quad (\mod 2).
\end{equation}
}
\end{itemize}

Indeed, since the Boole-coloring of each vertex is the sum of the
Boole-colorings of its incident edges, and recalling again Equation
(\ref{encaixada}), we can write
$$
\sum_{v\in V}\Psi (v) =
\sum^{3}_{i=1}n_{i}\1_{i}+(n-n')\0=\sum^{3}_{i=1}n_{i}\1_{i}=\sum^{}_{e\in
E}2\Psi (e) = \0,
$$
but this equality is only satisfied if $n_{i}\1_{i}=\0$ or
$n_{i}\1_{i}=\1_{i}$, for every $i\in \cal{C}$. Then, from
$n_{1}+n_{2}+n_{3}=n'$, we get the result.

Note that, as a direct consequence, we also get the following:
\begin{itemize}
  \item {\it There is no edge-coloring of a graph $G$ having only one vertex with
Boole-coloring $\1$ $($and the other vertices with Boole-coloring
$\0)$}.
\end{itemize}
Another consequence is the following result by Tait~\cite{Ta1880}:
\begin{itemize}
  \item {\it A cubic planar graph is Tait-colorable if and only if its corresponding map is $4$-colorable}.
\end{itemize}

\begin{figure}[t]
\begin{center}
\includegraphics[width=11cm]{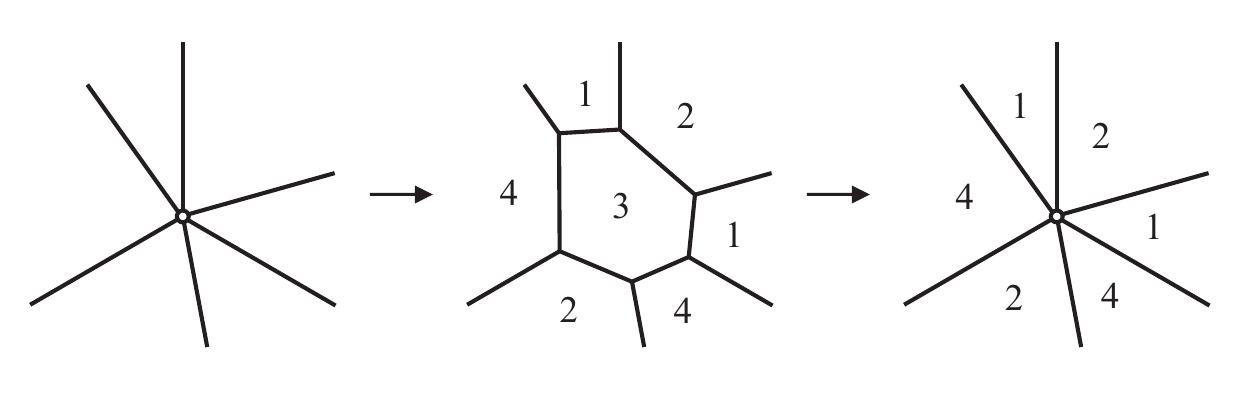}
\caption{An example of the fact that every map has a 3-graph
associated.} \label{fig:mapa=3graf}
\end{center}
\end{figure}

\begin{figure}[t]
\begin{center}
\includegraphics[width=8cm]{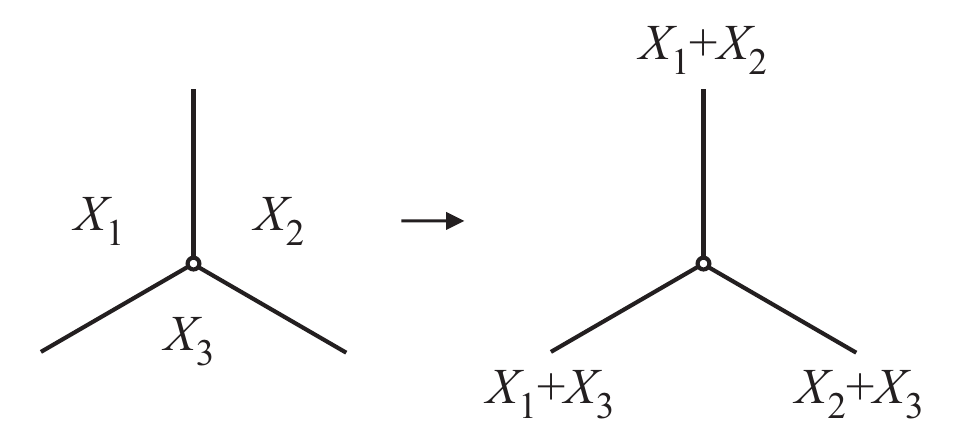}
\caption{Obtaining a Tait-coloring of a 3-graph.}
\label{fig:Tait-coloracio}
\end{center}
\end{figure}

\begin{figure}[t]
\begin{center}
\includegraphics[width=9cm]{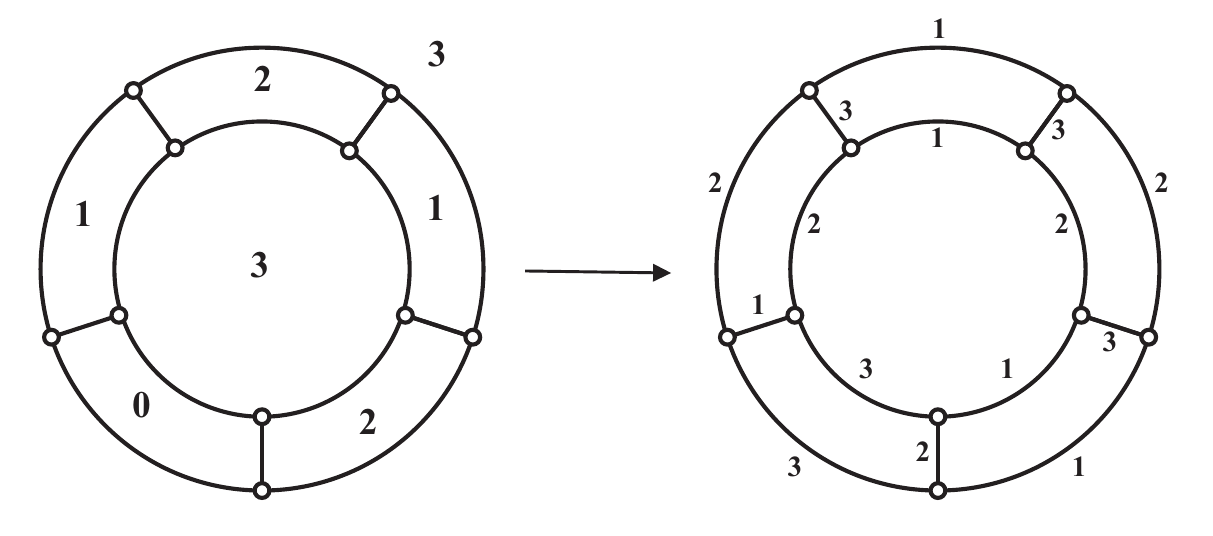}
\caption{The 4-coloring of a map and the Tait-coloring of its
edges.} \label{fig:colors_dodecaedre_roda}
\end{center}
\end{figure}

Using the Boole-colorings, the proof of this last result is as follows: First, recall that every map has a 3-graph
associated, because a vertex with degree greater than 3 can be replaced by a
polygon, in such a way that the map obtained can be colored with 4 colors, and so can the original
map; see an example in Figure~\ref{fig:mapa=3graf}.
Now assume that we have the regions of the map with the colorings
$\vec0,\1_1,\1_2,\1_3$. Then, to obtain a Tait-coloring of a cubic
planar graph, we only need to assign to each edge the sum of the
colorings of both regions separated by this edge. To see that this
gives a Tait-coloring, we only have to study one vertex, as
shown in Figure~\ref{fig:Tait-coloracio}. Since we have a 4-colored map, each two neighboring regions have different colors. Thus, no sum can give  $\vec0$. Moreover, since the
three regions with a common vertex have different colorings
$X_1,X_2$ and $X_3$ and $(\cal{B},+)$ is a group, the colorings
$X_1+X_2$, $X_1+X_3$ and $X_2+X_3$ must also be different.
Figure~\ref{fig:colors_dodecaedre_roda} provides an example of a
4-coloring of a map and its Tait-coloring (obtained from Table~1), where the colorings
$\vec0, \1_1, \1_2$ and $\1_3$ are denoted by $0,1,2$ and $3$,
respectively.

\begin{figure}[t]
\begin{center}
\includegraphics[width=13cm]{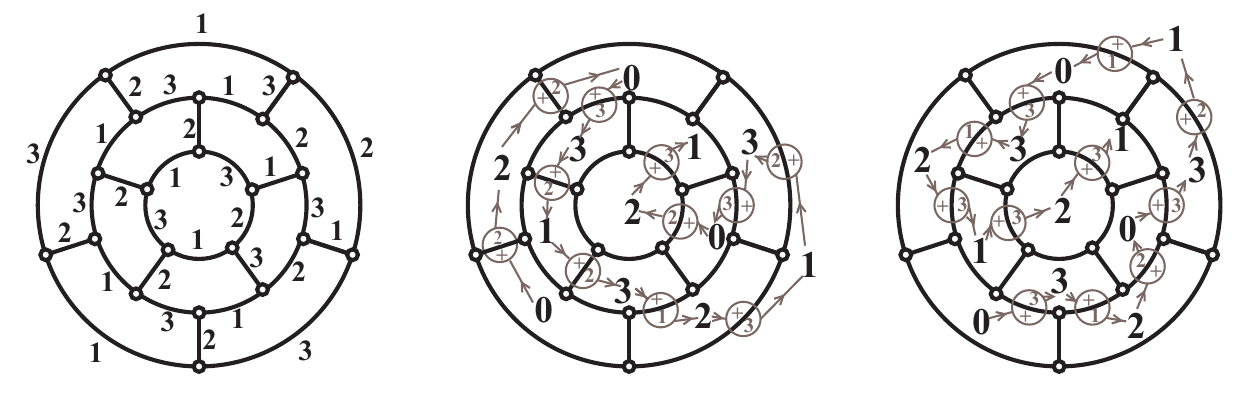}
\caption{An edge-coloring of the dodecahedron (also in Figure~\ref{fig:solids_platonics}) and two paths with the same
initial and final regions.} \label{fig:dodecaedre_roda}
\end{center}
\end{figure}

\begin{figure}[t]
\begin{center}
\includegraphics[width=5cm]{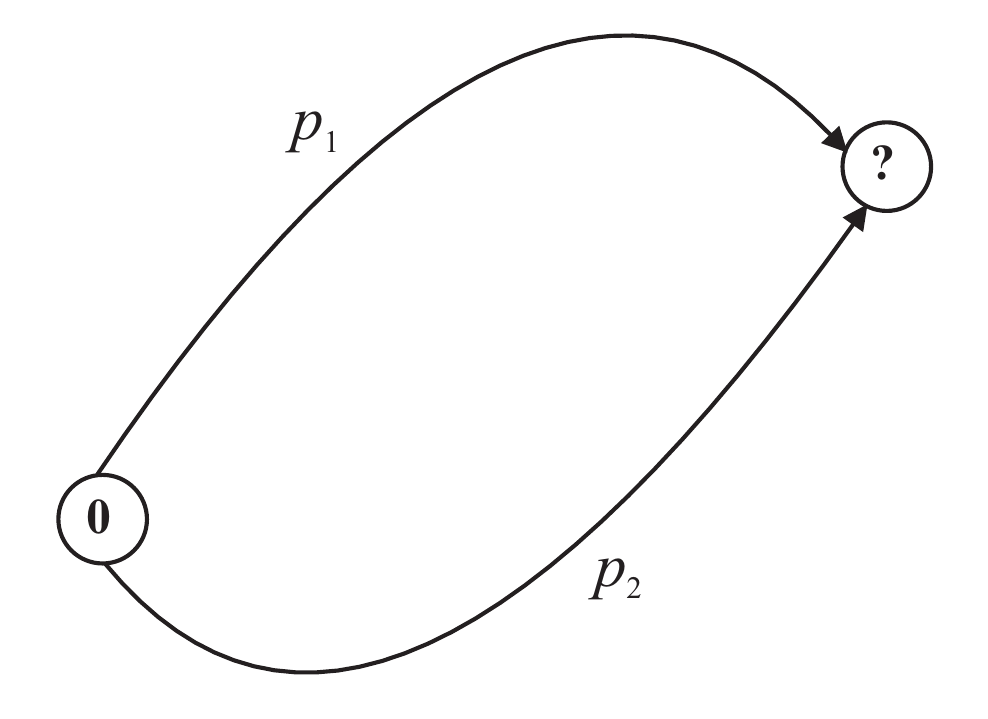}
\caption{Two paths from a region $\vec0$ to another with an unknown
color.} \label{fig:2camins}
\end{center}
\end{figure}

Conversely, if we want to obtain a
4-colored map from a Tait-coloring of the edges of the corresponding
graph, we begin by giving the coloring $\vec0$ to any region considered as initial.
Then, starting from this region, we follow an arbitrary path crossing some edges and visiting all the regions.
We give each newly visited region the coloring obtained by adding the coloring of the `previous' region plus the coloring of the last edge crossed. As no edge has the coloring $\vec0$, it is obvious that the coloring obtained for each region is different from that of its `previous' region in the path followed; for an example of this process, see Figure~\ref{fig:dodecaedre_roda} (left and center). Now, to
finish the proof, we need to show that the coloring of each region is
independent of the path followed. With this aim, let $p_1$ and $p_2$
be two paths with the same initial and final regions. We want to
prove that the coloring obtained for the final region is the same
following both paths; there is an example of this fact in
Figure~\ref{fig:dodecaedre_roda} (center and right). The colorings $X$ and $Y$ obtained by following both paths are equal if and only if the sum of the colorings of all edges crossed, respectively, by $p_1$ and $p_2$ is $\vec0$. Indeed, let $X_1,X_2,\ldots,X_s$ and $Y_1,Y_2,\ldots,Y_t$ be the colorings of the edges crossed respectively by $p_1$ and $p_2$, then $X_1+X_2+\cdots+X_s=X$ and $Y_1+Y_2+\cdots+Y_t=Y$. If $(X_1+X_2+\cdots+X_s)+(Y_1+Y_2+\cdots+Y_t)=\vec0$, the sums in both
parenthesis are equal, so $X=Y$. To prove this equality, we can
assume that $p_1+p_2$ is a simple curve (see
Figure~\ref{fig:2camins}) because, otherwise, we could
decompose it into some simple curves. If we imagine that we cut the
graph with this curve, we obtain two graphs, such that the colorings of the edges crossed by the curve must satisfy
$m_1 \equiv m_2 \equiv m_3 \ (\mod 2)$, where $m_i$ is the number of edges crossed with coloring $\1_i$. (Just imagine that in every cut we have two vertices of degree 1 and apply (\ref{parity-lemma}).) Then,
$(X_1+X_2+\cdots+X_s)+(Y_1+Y_2+\cdots+Y_t)=
m_1\,\1_1+m_2\,\1_2+m_3\,\1_3=\vec0$,
as claimed.

\begin{figure}[t]
\begin{center}
\includegraphics[width=4cm]{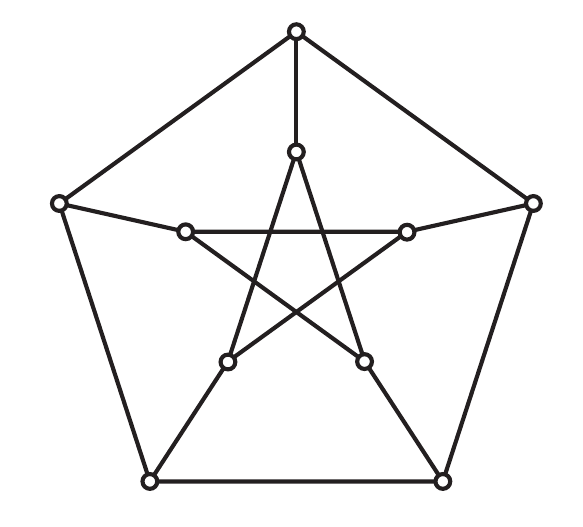}
\caption{The Petersen graph $P$.} \label{fig:petersen}
\end{center}
\end{figure}

\begin{figure}[t]
\begin{center}
\includegraphics[width=13.5cm]{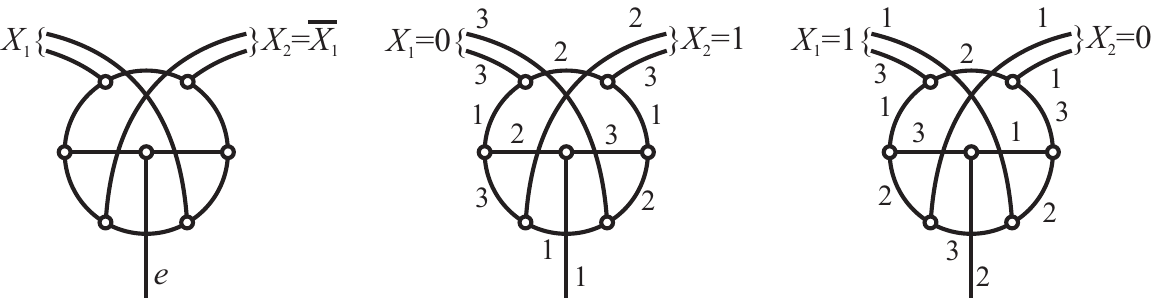}
\caption{Multipoles and the NOT gate.}
\label{fig:multipols}
\end{center}
\end{figure}

\begin{figure}[t]
\begin{center}
\includegraphics[width=5cm]{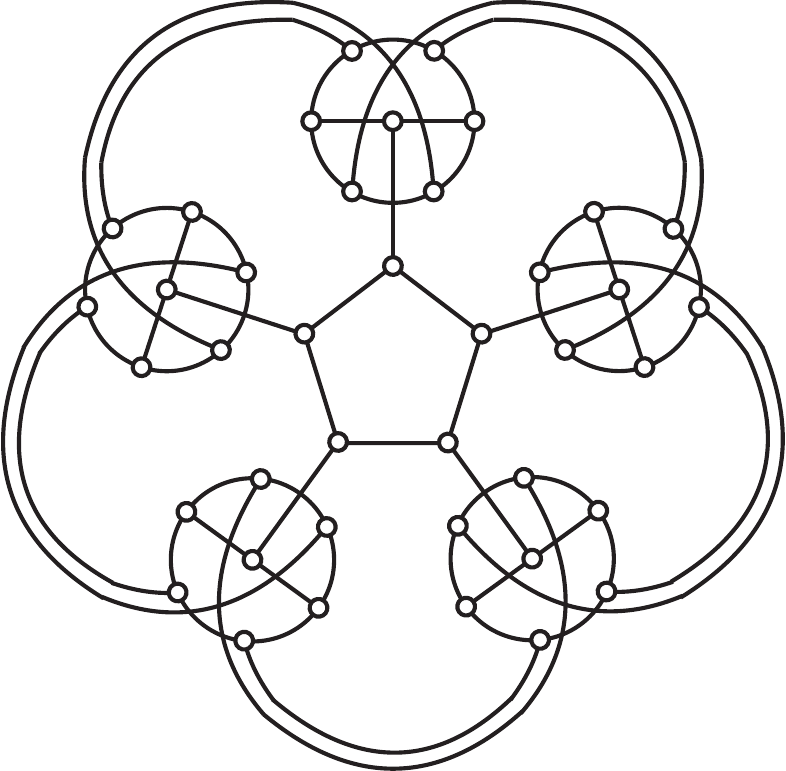}
\caption{A flower snark.}
\label{fig:flower-snark}
\end{center}
\end{figure}

As previously mentioned, the concept of colorings allows us to use the theory of Boolean algebra for the construction and characterization of \emph{snarks}, that is, cubic graphs that are not Tait-colorable, also known as \emph{class two}. The name `snark' was proposed by Gardner~\cite{Ga76}, who borrowed it
from a nonsense poem by the famous English author Lewis Carroll~\cite{Ca74}. The most simple example of snark is the
Petersen graph~\cite{Pe1898} (see Figure~\ref{fig:petersen}). With the colorings we can obtain infinite families of snarks. An example
is the family obtained by joining adequately an odd number of copies
of the multipole (cubic graph with edges and semi-edges---or `dangling edges'--- which are edges with only one final vertex), shown in Figure~\ref{fig:multipols} (left).
This structure behaves as a NOT gate of logic circuits
in the sense that, its edges and semi-edges having been Tait-colored,
the colorings $X_1$ and $X_2$ are conjugated one to each other, namely $X_2=\vec0$ (respectively, $X_2=\1$) if and only if $X_1=\1$ (respectively, $X_1=\vec0$). This is satisfied for any coloring of semi-edge $e$. Two examples of this fact are shown in  Figure~\ref{fig:multipols} (center and right). If, as previously stated,
we join an odd number of these multipoles in a circular configuration, adding some vertices to connect semi-edges $e$, any attempt at Tait-coloring
will lead to a conflict, and hence the graph is a snark. An example with five multipoles can be seen in Figure~\ref{fig:flower-snark}.
This family of snarks, called \emph{flower snarks}, was proposed by
Loupekhine (see Isaacs~\cite{i76}). The first infinite families
of snarks were given by Isaacs~\cite{i75}, but they can also be obtained by using Boole-colorings. More details on this technique can
be found in Fiol~\cite{f91}.


\section{Known and unknown: Ramsey theory}

Let us consider the following result:
\begin{itemize}
  \item {\it At a cocktail party with six o more people, there are always three people who are known or unknown to each other}.
\end{itemize}

In other words, if the complete graph $K_n$ on $n\ge 6$ vertices can
be (free) edge-colored with two colors, say blue and red, then it always contains a monochromatic triangle, namely, a subgraph $K_3$ with its three edges blue or red. Indeed, as each vertex $u$ has degree 5, at least 3 of its incident edges $\{u,v_i\}$, $1\le i\le 3$, must have the same color, for example, blue. Then, if any of the 3 edges $\{v_i,v_j\}$ ($1\le i<j\le 3$) is blue, we obtain a blue triangle. Otherwise, we have a red triangle.
Although this is an easy proof, it can be extremely difficult to prove similar results having more colors and/or imposing other monochromatic subgraphs. In this context, recall that, given $m$ graphs
$G_1,G_2,\ldots,G_m$, the {\it Ramsey number} $R(G_1,G_2,\ldots,G_m)$ is defined as the smallest number $n$, such that, in any edge-coloring of $K_n$ using $m$ colors, there always exists a monochromatic subgraph (with color $i$) isomorphic to $G_i$ for some $1\leq i\leq m$. If $G_i$ is a complete graph $K_r$, the Ramsey number is expressed by writing $r$ instead of $K_r$, for sake of simplicity. Some known results of exact values and bounds for Ramsey numbers are the following:
\begin{eqnarray*}
&& \hskip-.5cm R(3,3)=6, \ R(3,4)=9, \ R(3,5)=14, \ R(3,6)=18, \ R(4,4)=18,\\
&& \hskip-.5cm R(4,5)=25, \ 43\leq R(5,5)\leq 49; \ R(3,3,3)=17; \ 51\le R(3,3,3,3) \le 62.
\end{eqnarray*}
So, the result at the beginning of this section can be expressed as
$R(3,3)$ $\leq6$. Moreover, since $R(3,3)\geq 6$ (it is easy to color with two colors the edges of the complete graph
$K_5$ without monochromatic triangles: the `outer cycle' with one color and the 'inner' cycle with the other) we conclude that $R(3,3)=6$.
A good updated summary on this subject can be found in Radziszowski~\cite{r06}.

As an example, we now prove the following result:
\begin{itemize}
  \item
$R(3,3,3)=17$.
\end{itemize}
We first see that $R(3,3,3)\leq17$. We make an edge-coloring of a complete graph using three colors; say blue, red and green. Let us assume that the edge-coloring has no monochromatic triangles. The green neighborhood of a vertex $v$ is the set of vertices that have a green edge to $v$. The green neighborhood of $v$ cannot contain any green edge in order to avoid monochromatic triangles. Then, the edge-coloring of the green neighborhood of $v$ has only two colors: blue and red. Since $R(3,3)=6$, the green neighborhood of $v$ can contain at most 5 vertices. With the same reasoning, the blue and the red neighborhoods of $v$ can have at most 5 vertices each. As every vertex different from $v$ is in the green, blue or red neighborhoods of $v$, then the complete graph can have at most $1+5+5+5=16$ vertices. Thus, $R(3,3,3)\leq17$.

\begin{figure}[t]
\begin{center}
\includegraphics[width=12cm]{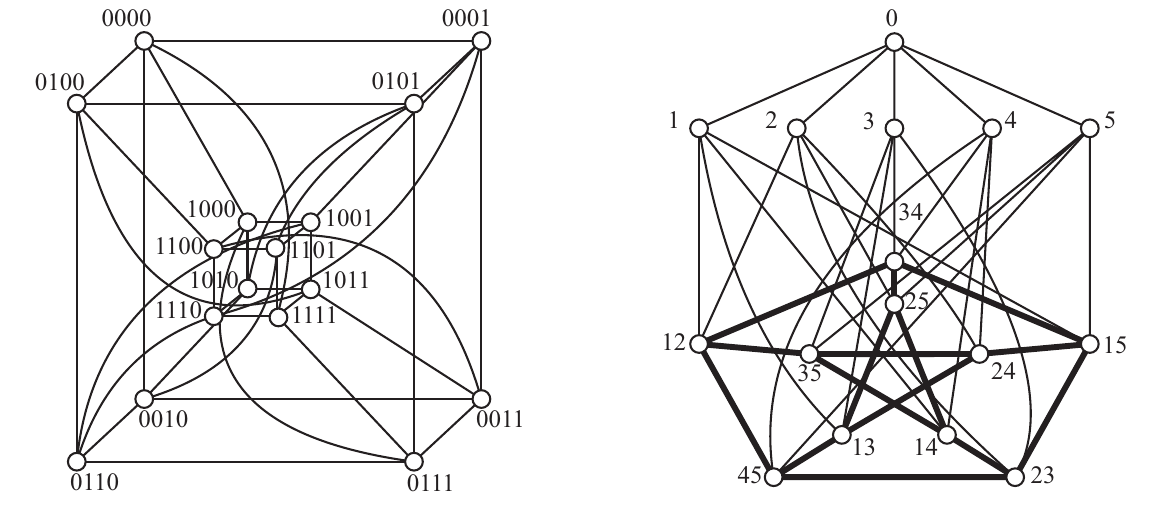}
\caption{Clebsch graph defined in two different ways.}
\label{fig:clebsch}
\end{center}
\end{figure}

Now, to prove that $R(3,3,3)\geq17$, we use algebraic graph theory based on the properties of eigenvalues and eigenvectors of the \emph{adjacency matrix}, that is, a matrix with rows and columns indexed by the vertices of the graph, and whose entries are either 1 or 0, according to whether the corresponding vertices are adjacent or not.

A $\delta$-regular graph with $n$ vertices is said to be \emph{$(n,\delta; a,c)$-strongly regular} if each pair of adjacent vertices has $a$ common neighbors and each pair of nonadjacent vertices has $c$ common neighbors.

If $R(3,3,3)\geq17$, then we can color the edges of the complete graph $K_{16}$ with three colors, namely, we can make an edge-coloring of $K_{16}$ without monochromatic triangles. The required edge-coloring is equivalent to a decomposition of $K_{16}$ into three graphs $G_1$, $G_2$ and $G_3$, each one corresponding to one color. It follows that each $G_i$, $i=1,2,3$, must be a graph on 16 vertices, regular of degree 5 (because each vertex has degree 15 and the neighborhood with one color has at most 5 vertices) and without triangles. Moreover,  each vertex $u\in V_i$ has 10 vertices at distance 2, which can be reached by $5\cdot 4=20$ paths of length 2. Then, we can consider a graph in which any two nonadjacent vertices have 2 common neighbors
and any two adjacent vertices have no common neighbors. In other words, a $(16,5;0,2)$-strongly regular graph. It is known that there is just one such graph, the Clebsch graph, which is illustrated in two
different ways in Figure~\ref{fig:clebsch}. On the left, there is
the Clesbch graph, as the graph whose vertices are labeled with the
numbers 0 to 15 in base 2, and where two vertices are adjacent whenever
the corresponding labels differ either by one or by all four digits.
On the right, there is the Clebsch graph, as the rooted graph with
vertices labeled $0,i$, and the unordered pairs $ij$, with
$i,j\in\{1,2,3,4,5\}$, for $i\not=j$. In this representation, the
adjacencies are $0\sim i$, $ij\sim i$, $ij\sim j$, and $ij\sim
kl$ if $i,j,k,l$ are all different and $i,j,k,l\in\{1,2,3,4,5\}$. In
fact, the Clebsch graph is vertex-transitive (informally speaking, we see the same structure from any vertex), so that any vertex can be chosen as vertex 0. Notice that, from this view of the Clebsch graph, it is apparent that the induced subgraph on ten vertices at
distance 2 (from the vertex chosen as 0) is the Petersen graph~\cite{Pe1898}; compare Figure~\ref{fig:clebsch} (on the right) and Figure~\ref{fig:petersen}.

Therefore, our problem is to find three edge-disjoint copies of the
Clebsch graph in $K_{16}$. To this end, let us introduce the following terminology: Let $G_i=(V,E_i)$ be a family of graphs on the same
vertex set $V$ and such that $E_i\cap E_j=\emptyset$, for
$i,j=1,2,\ldots,m$. We define the graph $G=\bigcup_{i=1}^m G_i$ as
the graph $G=(V,E)$, where $E=\bigcup_{i=1}^m E_i$. Notice that the
corresponding adjacency matrices satisfy $\A(G)=\sum_{i=1}^m \A(G_i)$. With $Cl_i$ denoting a graph isomorphic to the Clebsch graph, our problem now reads: Is it true that $K_{16}=Cl_1\cup Cl_2\cup Cl_3$? In terms of their adjacency matrices $\A_i=\A(Cl_i)$, we have
\begin{equation}
\label{j-i}
\A_1+\A_2+\A_3=\J-\I,
\end{equation}
since the adjacency matrix of $K_{16}$ is equal to $\J-\I$, where $\J$ denotes the matrix whose entries are all 1 and $\I$ is the identity matrix.

We now use eigenvalue techniques to address Equation (\ref{j-i}).
Recall that the spectrum of an adjacency matrix gives the
eigenvalues of this matrix (which are real because the matrix is symmetric), and that each eigenvalue has at least one eigenvector associated. To find the spectra of the Clebsch graph
and the matrix $\J-\I$, we can either compute them or simply find them in some standard reference, such as Godsil and Royle~\cite{GoRo01}. We then have that $\sp \A_i=\{5^1,1^{10},-3^5\}$ and $\sp (\J-\I)=\{15^1,-1^{15}\}$, where the superscripts denote the multiplicity of each eigenvalue. In both cases, the largest eigenvalue has the all-1 vector $\j$ as eigenvector. It follows that the eigenvectors of the other eigenvalues are in the subspace ${\cal H}=\j^\bot$ (with vectors the addition of whose components are zero). Denote by ${\cal E}_i$ the eigenspace of  $\A_i$ corresponding to the eigenvalue 1, namely, ${\cal E}_i=\ker (\A_i-\I)$, and consider the subspace ${\cal F}={\cal E}_1\cap{\cal E}_2\subset {\cal H}$. As $\dim {\cal E}_1=\dim {\cal E}_2=10$ and $\dim {\cal H}=15$, we infer that $\dim {\cal F}\geq 5$. From Equation (\ref{j-i}), with $\A_1\v=\v$, $\A_2\v=\v$ and $(\J-\I)\v=-\v$, where $\v\in {\cal F}$, we obtain that $\A_3\v=-3\v$ and, then, $\dim {\cal F}=5$ and ${\cal F}=\ker(\A_3+3\I)$. This implies that
$$
{\cal H}={\cal F}_1\cup {\cal F}_2\cup {\cal F}_3
$$
where ${\cal F}_i={\cal E}_j\cap{\cal E}_k$, with $\{i,j,k\}=\{1,2,3\}$.

\begin{figure}[t]
\begin{center}
\includegraphics[width=6cm]{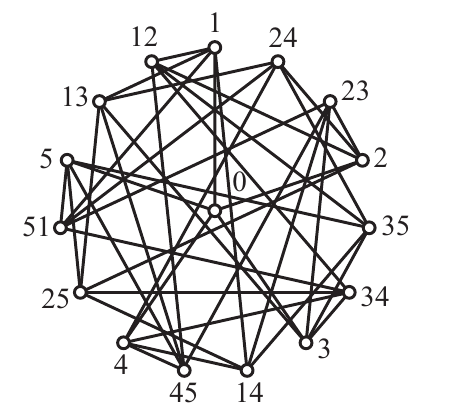}
\caption{$K_{16}/3$ $=$ Clebsch graph.}
\label{fig:clebsch2}
\end{center}
\end{figure}

This indicates that the required spectral condition necessary to the existence of the decomposition $K_{16}=Cl_1\cup Cl_2\cup Cl_3$ is satisfied. In this case, this condition is also sufficient, and it is known that there are only two nonisomorphic decompositions. One of these is illustrated in Figure~\ref{fig:clebsch2}, which shows how to color one third of the edges of $K_{16}$ with one color using
the Clebsch graph. By rotating this graph $\frac{2\pi}{15}$ and $\frac{4\pi}{15}$ radians, we obtain the edges to be colored with the two other colors; with this, we get $R(3,3,3)=17$.

In the case of avoiding monochromatic triangles with $m>3$ colors, only bounds of Ramsey numbers are known. By definition, we state that  $C(m):=R(3,3,\stackrel{m}{\ldots},3)-1$ for $m\geq1$, that is, $C(m)$ is the biggest integer $n$ such that $K_n$ can be colored with $m$ colors without monochromatic triangles. The following upper bound is known (see Fiol, Garriga and Yebra~\cite{fgy95}):
\vskip -.2cm
\begin{equation}
\label{fita}
\hskip -8.8cm \bullet\ \ C(m)\le \lfloor m! \, e\rfloor,
\end{equation}

Recall that, surprisingly, we find the number $e$. The proof is as
follows: Obviously, $C(1)=R(3)-1=2$ and we know that $C(2)=R(3,3)-1=5$ and $C(3)=R(3,3,3)-1=16$. If we compute $C(3)$ from $C(2)$, considering that a vertex $v$ can only be adjacent to $6+5+5$ vertices, we obtain that $C(3)\leq 3C(2)+1=16$. For any $m\geq1$, we get the recurrence
$$
C(m+1)\leq (m+1)\,C(m)+1.
$$
We solve the corresponding linear equation
$$
D(m+1) = (m+1)\,D(m)+1,
$$
first solving its homogeneous equation
$$
D(m+1) = (m+1)\,D(m) \Rightarrow D(m) = K\,m!,
$$
where $K$ is a constant. Then, we look for a particular solution
$D(m)=K(m)\,m!$ of the complete equation:
\begin{eqnarray*}
&& \hskip-.5cm K(m+1)(m+1)!=(m+1)K(m)m!+1\\
&& \hskip-.5cm \Rightarrow\quad K(m+1)-K(m)=\frac{1}{(m+1)!}
\quad\Rightarrow\quad K(m)=\sum_{r=1}^{m} \frac{1}{r!}+\alpha\\
&& \hskip-.5cm \Rightarrow\quad D(m)=m!\left(\sum_{r=1}^{m}
\frac{1}{r!}+\alpha\right),
\end{eqnarray*}
where $\alpha$ is a constant. Finally, $C(1)=D(1)=2$ gives
$\alpha=1$ and, hence, $C(m)\leq \lfloor m!\,e\rfloor$, as claimed.

From the examples given at the beginning of this section, we saw that
$51\le R(3,3,3,3) \le 62$. Using (\ref{fita}), we obtain that
$$
R(3,3,3,3)=C(4)+1\leq \lfloor 4!\,e\rfloor +1 = 66,
$$
which represents a good upper bound, quite close to the best bound known.

\section{Common friends: Distance-regularity and coding theory}

As commented by Aigner and Ziegler~\cite{az98}, nobody knows who was
the first to state the following result and to give it the human
touch:
\begin{itemize}
\item {\it At a cocktail party with three or more people, if each two people have exactly one friend in common, then there is a person $($the `politician'$)$ who is a friend of everybody}.
\end{itemize}

\begin{figure}[t]
\begin{center}
\includegraphics[width=7cm]{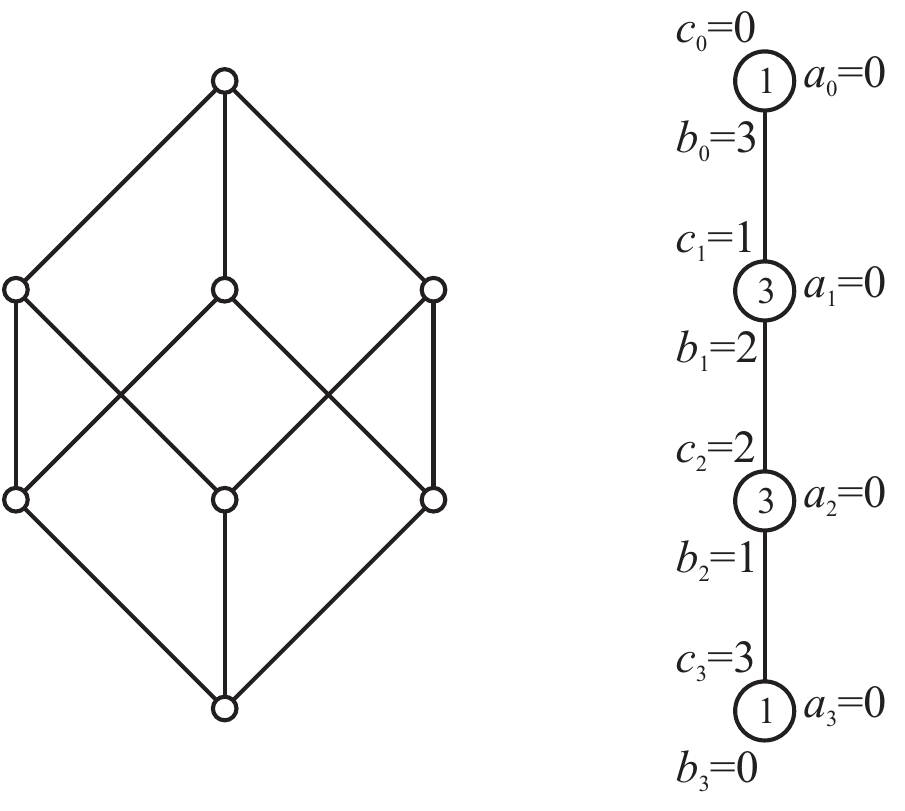}
\caption{A layer partition of the cube $Q$ and its intersection
diagram.} \label{fig:cub_capes}
\end{center}
\end{figure}

Nowadays, this result is known as the \emph{Friendship theorem}. As
mentioned in the introduction, the first proof (by contradiction)
was given by Erd\H{o}s, R\'enyi and S\'os~\cite{ers66} in 1966, and
is considered to be the most successful. Basically, it has two
parts: First, it is proved that if the graph $G$ which models such
a cocktail party (where people correspond to vertices and  friendships are represented by edges) is a
counterexample with more than three vertices, then it has to be
regular, say with degree $k$. As a consequence, $G$ has to be
\emph{strongly regular} with parameters $(n,k;1,1)$, that is, every
two adjacent vertices has exactly one common neighbor, and the same
holds for every two nonadjacent vertices. Second, spectral graph theory is used to prove that $G$ cannot exist. In fact, the
hypothetic graph $G$ would be an example of a distance-regular
graph, in this case with diameter 2 (the concepts of strongly-regularity and distance-regularity coincide for connected graphs with diameter 2). Generally speaking, we say that
a graph is \emph{distance-regular} if, when it is observed or `hung'
from any of its vertices (called \emph{root}), we obtain a partition of the vertex set into layers, where the layer $i$ contains the vertices at distance $i$ from the root, and the vertices in a layer are indistinguishable from each other with respect to their adjacencies. A more precise definition of distance-regularity is the following: A graph $G$ with diameter $D$ is distance-regular if, for every pair of vertices $u,v$ and integers $0\le i,j\le D$, the number $p_{ij}(u,v)$ of vertices at distance $i$ from $u$ and at distance $j$ from $v$ only depends on the distance between $u$ and $v$, $\dist(u,v)=k$. Then, we write $p_{ij}(u,v)=p_{ij}^k$, where the constants $p_{ij}^k$ are called the \emph{intersection numbers}.
Indeed, because of the many relations between these numbers, it is
possible to give a much more simple definition, since for each
distance $k$ we only need the pairs of distances $(i,j)=(k-1,1)$,
$(k,1)$ and $(k+1,1)$. The corresponding intersection numbers are
enough to determine all the others; see, for example,
Biggs~\cite{b93}. Therefore, the most common definition of
distance-regularity is: A graph $G$ is distance-regular if, for
every pair of vertices $u,v$ at distance $\dist(u,v)=k$, the numbers
$c_k$, $a_k$, and $b_k$ of vertices adjacent to $v$, and at distance
$k-1$, $k$, and $k+1$, respectively, from $u$ only depends on $k$,
such that $c_k=p_{k-1,1}^k$, $a_k=p_{k,1}^k$, and $b_k=p_{k+1,1}^k$.
As simple examples of distance-regular graphs, we have the
$1$-skeleton of regular polyhedrons; see again
Figure~\ref{fig:solids_platonics}. In Figure~\ref{fig:cub_capes}, we
show the layer partition of the cube graph $Q$ with the so-called \emph{intersection diagram} of the corresponding intersection numbers. Notice that each layer is represented by a circle containing its number of vertices.

Since their introduction by Biggs in the early 70's,
distance-regular graphs, and their principal generalization called
\emph{association schemes} (see, for example, Brouwer and Haemers~\cite{BrHa95}), have been key concepts in algebraic
combinatorics. These graphs have connections with other
areas of mathematics, such as geometry, coding theory, group theory,
design theory, and with other parts of graph theory. As pointed out by Brouwer, Cohen and Neumaier in their monumental book on this
subject~\cite{bcn89}, this is because most of the finite objects
with `enough' regularity are closely related to distance-regular graphs.

In 1997 Fiol and Garriga~\cite{fg97,f99} gave the following quasi-spectral characterization of distance-regular graphs:
\begin{itemize}
\item {\it A regular graph $G$ with adjacency matrix $\A$ and $d+1$ distinct eigenvalues is distance-regular if and only if the number $|\Gamma_d(u)|$ of vertices at distance $d$ from each vertex $u$ is a constant and only depends on the spectrum of the matrix $\A$}.
\end{itemize}
More precisely, consider a regular graph $G$ with $n$ vertices and spectrum $\sp G = \{\lambda_0^1,\lambda_1^{m_1},\ldots,\lambda_d^{m_d}\}$,
where $\lambda_0,\lambda_1,\ldots,\lambda_d$ are the eigenvalues of
$\A$ and the superscripts denote their multiplicities; $\lambda_0$ is
simple because $G$ is connected, thus $\A$ is irreducible (Perron-Frobenius theorem for nonnegative matrices, see Godsil~\cite[p. 31]{g93}). Then, $G$ is distance-regular if and only if, for each vertex $u$,
\begin{equation}
\label{numcharac} |\Gamma_{d}(u)| =
n\left(\sum_{i=0}^{d}\frac{\pi_0^2}{m_i \pi_i^2}\right) ^{-1},
\end{equation}
where $\pi_i$'s are moment-like parameters, which can be calculated
from the distance between eigenvalues with the formula
$\pi_i=\prod_{j=0(j\neq i)}^{d}|\lambda_i-\lambda_j|$, for $0\le
i\le d$. As examples, we give the spectrum, the number of vertices
and the value of $|\Gamma_{d}(u)|$ obtained from Equation (\ref{numcharac}) of the cube $Q$ and the Petersen graph $P$ (see again Figures~\ref{fig:cub_capes} and~\ref{fig:petersen}, respectively):
\begin{itemize}
  \item Cube: $\sp Q=\{3^1, 1^3, -1^3, 3^1\}$, $n=8$, $|\Gamma_3(u)|=1$.
  \item Petersen: $\sp P=\{3^1,1^5,-2^4\}$, $n=10$,
  $|\Gamma_2(u)|=6$.
\end{itemize}

As previously mentioned, the theory on distance-regular graphs
has many applications in coding theory. Recall that a code $C$, with a set of allowed words or {\it code-words}, can be simply represented as a vertex subset of a distance-regular graph $G$; see Godsil~\cite{g93} and van Lint~\cite{vl99}. The vertex subset represents the `universe' of words, with or without meaning,
which can be received. There is an edge between two words
if, with a certain probability, one can be transformed into the
other in the process of transmission. Then, the shorter the distance
between two words in $G$, the more similar the words. If a
code-word has not suffered too many changes, the resulting word is
not far from the original one and it is possible to retrieve it
(decision criterion by proximity). Therefore, a code is better if
the words that constitute it are far away from each other. In the
study and design of good codes, some algebraic techniques are used
to obtain information about the structure of the graph $G$ and, in particular, about the vertex subset $C$ that represents the code. In the applications of special relevance, there are the so-called \emph{completely regular codes}, whose graphs are structured in a kind of distance-regularity around the set that constitutes the code. Thus, these codes can be algebraically characterized in a similar way to the characterization of the distance-regular graphs through their spectra; see Fiol and Garriga~\cite{fg99} for more information.

\section{Weddings: Hall's and Menger's theorems. Multibus networks}

Let us imagine two groups of heterosexual people available for marriage, one of women and another of men, the latter at least as large as the former. Also imagine that every woman knows a certain number of men. The Hall Marriage theorem gives necessary and sufficient conditions for every woman to be able to marry a man who she knows:
\begin{itemize}
\item {\it A \emph{complete matching} is possible if and only if each group of women, whatever their number, knows altogether at least an equal number of men}.
\end{itemize}
If the sets of women and men are denoted by $U$ and $V$, respectively,  we can represent the above situation as a bipartite graph $G=G(U\cup V,E)$, with stable vertex sets $U$ and $V$ and where edges stand for acquaintances. Then, we can state Hall's theorem in a more mathematical form:
\begin{itemize}
\item {\it In a bipartite graph $G=G(U\subset V,E)$ with $|U|\leq|V|$, a complete matching is possible if and only if, for every  $U^*\subset U$,
\begin{equation}
\label{Hall}
|\Gamma(U^*)|\geq|U^*|,
\end{equation}
where $\Gamma(U^*)=\cup_{u\in U^*}\Gamma(u)$}.
\end{itemize}
(Recall that $\Gamma(u)\subset V$ is the set of vertices adjacent to vertex $u\in U$.)

\begin{figure}[t]
\begin{center}
\includegraphics[width=8cm]{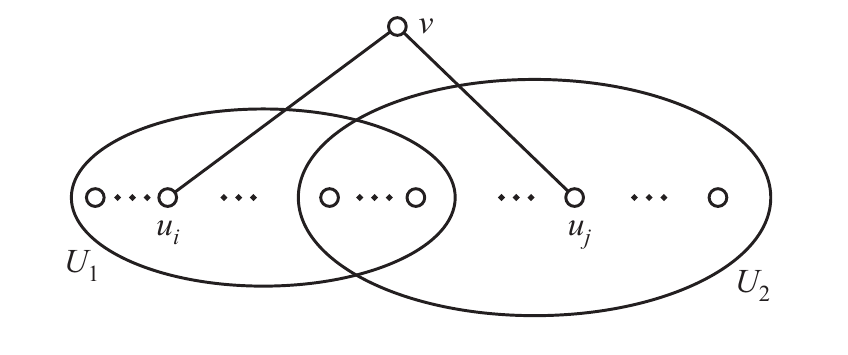}
\caption{The situation of the proof of Hall's theorem.} \label{fig:Hall}
\end{center}
\end{figure}

 There are several proofs of Hall's theorem. The proof we present here is by Rado, although our reasoning is a little different from that in Bollob\'{a}s~\cite{Bo04} or Harary~\cite{h69}. As necessity  is trivial, we are going to prove sufficiency. If graph $G$ satisfies Eq.~(\ref{Hall}), for any $u_i,u_j\in U$ with $i\neq j$ and $\Gamma(u_i)\cap\Gamma(u_j)=\emptyset$, it is immediate that $G$ contains a complete matching. If $\Gamma(u_i)\cap\Gamma(u_j)\neq\emptyset$, then there exist at least two edges $u_iv$ and $u_jv$, with $v\in V$. Now we claim that, after removing one of these edges, the resulting graph still satisfies Eq.~(\ref{Hall}). Indeed, if this were not the case, there would be two subsets $U_1,U_2\subset U$, with $u_i\in U_1$ and $u_j\in U_2$, such that $|\Gamma(U_1)|=|U_1|$ and $|\Gamma(U_2)|=|U_2|$. Moreover, $u_i$ would be the only vertex of $U_1$ adjacent to (some vertex of) $V$, and $u_j$ would be the only vertex of $U_2$ adjacent to $V$. See this situation in Figure~\ref{fig:Hall}. Then, we would have that the common number of adjacent vertices to $U_1$ and $U_2$ would satisfy the inequality:
\begin{eqnarray*}
|\Gamma(U_1)\cap\Gamma(U_2)|  & \geq  & |\Gamma(U_1-\{u_i\})\cap\Gamma(U_2-\{u_j\})|+1
   \geq  |\Gamma(U_1\cap U_2)|+1 \\
    & \geq & |U_1\cap U_2|+1.
\end{eqnarray*}
Moreover, we would also have:
\begin{eqnarray*}
  |\Gamma(U_1\cup U_2)| & = & |\Gamma(U_1)\cup\Gamma(U_2)| = |\Gamma(U_1)|+|\Gamma(U_2)|-|\Gamma(U_1)\cap\Gamma(U_2)| \\
  & \leq  & |\Gamma(U_1)|+|\Gamma(U_2)|-|U_1\cap U_2|-1 \\
   & = & |U_1|+|U_2|-|U_1\cap U_2|-1,
\end{eqnarray*}
a contradiction since,  according to Eq.~(\ref{Hall}),
$$
|\Gamma(U_1\cup U_2)| \geq |U_1\cup U_2| = |U_1|+|U_2|-|U_1\cap U_2|.
$$
Consequently, every vertex $v\in V$ with degree $\delta(v)\ge 2$ can be converted to a vertex with degree $1$, and the resulting graph still satisfies Eq.~(\ref{Hall}). This completes the proof.

Curiously, Hall's theorem is closely linked to another classical
result in graph theory: Menger's theorem; see, for example,
Bollob\'{a}s~\cite{Bo90}. As in the case of Hall's theorem, Menger's
theorem states that a certain condition, which is trivially
necessary for a result to be true, is also sufficient. In Menger's case,
the result is not on matchings, but on the vertex-connectivity $\kappa$ (or edge-connectivity $\lambda$) of a graph, which is defined as the minimum cardinality of a vertex (or edge) set whose deletion
disconnects the graph or, in particular, two given vertices $u,v$. This set is called a \emph{cutting set} or \emph{separating set} of $G$ or, in particular, of $u,v$. Then, Menger's theorem states that
for every pair of vertices $u,v$ (nonadjacent, in the case of
computing $\kappa$):
\begin{itemize}
\item \emph{The minimum size $\kappa(u,v)$ of a separating set of vertices equals the maximum number of independent paths in vertices from $u$ to $v$.}
\item \emph{The minimum size $\lambda(u,v)$ of a separating set of edges equals the maximum number of independent paths in edges from $u$ to $v$.}
\end{itemize}

It has been shown that the vertex-connectivity $\kappa=\min_{u,v\in V}\kappa(u,v)$ (or edge-connectivity $\lambda=\min_{u,v\in V}\lambda(u,v)$) of a graph or digraph $G$ (a digraph is a graph whose edges are
associated to one of the two possible directions) reaches its maximum value, which equals the minimum degree of $G$, if in $G$ the diameter is small enough with respect to the girth (see F\`{a}brega and Fiol~\cite{ff89}) or if the number of vertices is large enough with respect to the diameter (see Fiol~\cite{f93}).

\begin{figure}[t]
\begin{center}
\includegraphics[width=9cm]{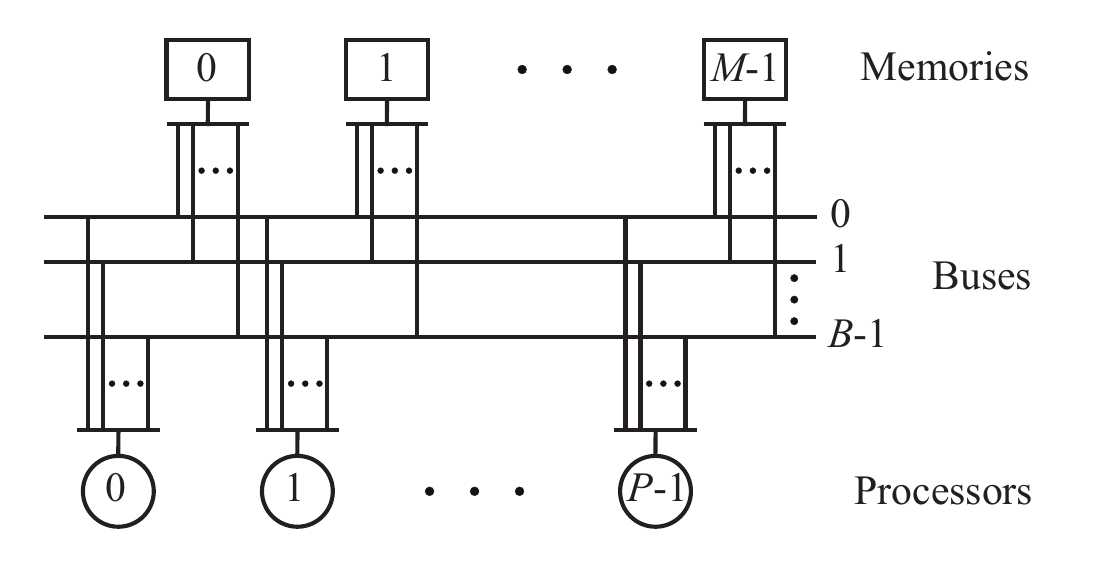}
\caption{The complete multibus interconnection scheme.}
\label{fig:multibus}
\end{center}
\end{figure}

Both the theorems mentioned, Hall's and Menger's, have many
applications in the study and design of interconnection networks
(for example, between processors) and in communication networks. Here
we explain an application of Hall's theorem to the study of
\emph{multibus interconnection networks}: A multiprocessor system
with shared memory and multibus interconnection network consists of
$P$ processors, $B$ buses and $M$ memory modules with
$B\le\min\{P,M\}$. The processors have access to the memory modules
through the buses, so we can establish processor-bus and bus-memory connections. Let us assume that there are $m\le
M$ requirements by the processors for accessing to different
memory modules. As each processor-memory connection
requires a bus, if $m\leq B$, then $m$ memories will be assigned;
instead, if $m>B$, then only $B$ memories will be assigned. In
the complete scheme (see Figure~\ref{fig:multibus}), each bus is
connected to all the memories and all the processors. This represents $B(P+M)$ connections, and generally this provides an important saving with respect to the \emph{crossbar} network with $PM$ connections, one connection between each pair processor-memory, because the number of buses is normally much smaller than the number of processors and memories. For example, if $M=N$ (an usual situation), the saving is obtained if $B<M/2$.

\begin{figure}[t]
\begin{center}
\includegraphics[width=10cm]{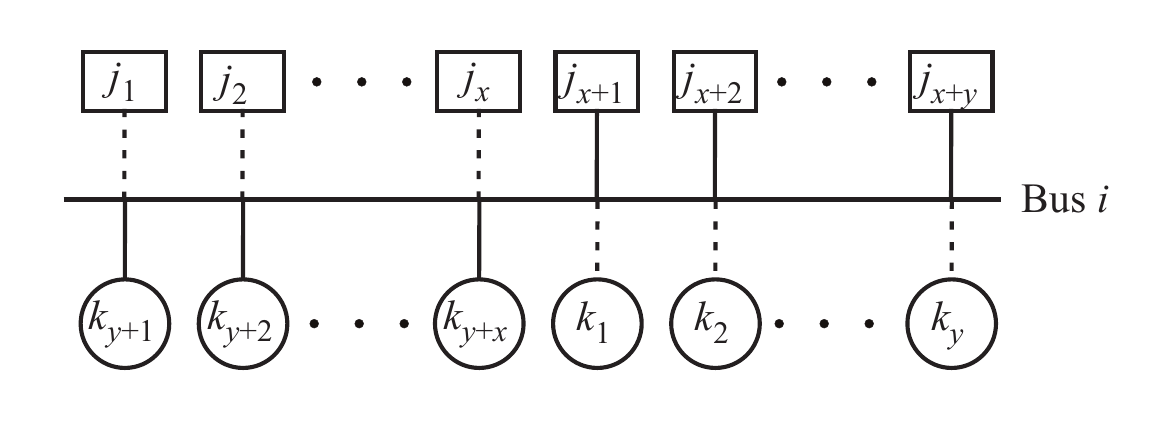}
\caption{Part of a system that suffers degradation.}
\label{fig:degradacio}
\end{center}
\end{figure}

Because the cost of the network basically depends on the number of connections, it is useful to consider the redundancy of this scheme. Namely, what is the maximum number of connections (processor-bus or bus-memory) that can be removed without having system degradation? In other words, {\it how many connections, from all of $B(P+M)$, can be removed such that any of the $m\le B$ processors asking for access to
any of the $m$ different memory modules do not lose access?} The answer is a direct consequence of the following result:

\begin{itemize}
\item
{\em In a multiprocessor system with multibus network without having degradation, each bus can be disconnected from at most $B-1$ altogether processors or memory modules.}
\end{itemize}
The proof is as follows: For each bus $i$, $0\leq i\leq B-1$, let $p_i$ and $m_i$ be, respectively, the number of processors and memories connected to it. Analogously, let $\overline{p}_i$ and  $\overline{m}_i$ be the numbers of processors and memories disconnected from bus $i$. Obviously, $p_i+\overline{p}_i=P$ and
$m_i+\overline{m}_i=M$. The result states that, in a non-degrading system, each bus $i$ can be disconnected from, at most, $B-1$ processors or memories, namely, $\overline{p}_i+\overline{m}_i\le B-1$ for $0\leq i\leq B-1$. But we can also state that each bus must have more than $P+M-B$ connections, such that $p_i+m_i>P+M-B$ for $0\leq i\leq B-1$. Assume that, on the contrary, for each bus $i$, we have $\overline{p}_i+\overline{m}_i\geq B$. Let $k_1,k_2,\ldots,k_y$
with $y\leq \overline{p}_i\leq P$ and $j_1,j_2,\ldots,j_x$ with $x\leq \overline{m}_i\leq M$ be, respectively, the processors and memories disconnected to the bus $i$. Note that $x+y=B$. Now consider $x$ other processors $k_{y+1},k_{y+2},\ldots,k_{y+x}$ and $y$ other memories $j_{x+1},j_{x+2},\ldots,j_{x+y}$, as in Figure~\ref{fig:degradacio}. Let $(k,j)$ be the requirement of  processor $k$ to access to memory $j$. None of the $B$ requirements
$$
(k_{1},j_{x+1}), (k_{2},j_{x+2}), \ldots, (k_{y},j_{x+y}), (k_{y+1},j_{1}), (k_{y+2},j_{2}), \ldots, (k_{y+x},j_{x})
$$
can use bus $i$, and this means that the system suffers degradation.

So, as stated before, the conclusion is that the maximum number of redundant connections is $B(B-1)$. In fact, this value is obtained
with the so-called \emph{minimum topologies}, such as the
\emph{rhombic} and the \emph{staircase} topologies; see Tables~2 and 3, respectively. More details can be found in Fiol, Valero, Yebra and Land~\cite{fvyl84} and in Lang, Valero and Fiol~\cite{lvf83}.

Notice that the result only gives us necessary conditions for suffering degradation. In this context, Hall's theorem is used to give a characterization for the interconnection topologies to prevent degradation of the system, as in the aforementioned cases of the complete and the minimum topologies:
\begin{itemize}
\item {\it A multibus system does not suffer degradation if and only if any of the $p\le B$ disjoint pairs processor-memory are connected to a set of, at least, $p$ buses}.
\end{itemize}

As previously stated, this result gives necessary and sufficient conditions for a non-degrading multibus system.

\begin{table}
\begin{center}
\label{taula:romb}
\begin{tabular}{|rrrrrrrrrrrrrrrr|} \hline
\multicolumn{16}{|c|}{Rhombic scheme}\\
\hline\hline
 0 &  1 & 2 & 3 & 4 & 5 & 6 & 7 & 8 & 9 & 10 & 11 & 12 & 13 & 14 & 15\\
\hline
  &   &   &   &   &   &   & 7 & 7 & 7 & 7 & 7 & 7 & 7 & 7 & 7\\
  &   &   &   &   &   & 6 & 6 & 6 & 6 & 6 & 6 & 6 & 6 & 6 &  \\
  &   &   &   &   & 5 & 5 & 5 & 5 & 5 & 5 & 5 & 5 & 5 &   &  \\
  &   &   &   & 4 & 4 & 4 & 4 & 4 & 4 & 4 & 4 & 4 &   &   &  \\
  &   &   & 3 & 3 & 3 & 3 & 3 & 3 & 3 & 3 & 3 &   &   &   &  \\
  &   & 2 & 2 & 2 & 2 & 2 & 2 & 2 & 2 & 2 &   &   &   &   &  \\
  & 1 & 1 & 1 & 1 & 1 & 1 & 1 & 1 & 1 &   &   &   &   &   &  \\
0 & 0 & 0 & 0 & 0 & 0 & 0 & 0 & 0 &   &   &   &   &   &   &  \\
\hline
\end{tabular}
\caption{Matrix representation of the rhombic scheme
with $M=16$ and $B=8$ (entries indicate the buses connected to
memory modules).}
\end{center}
\end{table}

\begin{table}
\begin{center}
\label{taula:escala}
\begin{tabular}{|rrrrrrrrrrrrrrrr|} \hline
\multicolumn{16}{|c|}{Staircase scheme}\\
\hline\hline
0 & 1 & 2 & 3 & 4 & 5 & 6 & 7 & 8 & 9 & 10 & 11 & 12 & 13 & 14 & 15\\
\hline
  &   &   &   &   &   &   & 7 & 7 & 7 & 7 & 7 & 7 & 7 & 7 & 7\\
  &   &   &   &   &   & 6 &   & 6 & 6 & 6 & 6 & 6 & 6 & 6 & 6\\
  &   &   &   &   & 5 &   &   & 5 & 5 & 5 & 5 & 5 & 5 & 5 & 5\\
  &   &   &   & 4 &   &   &   & 4 & 4 & 4 & 4 & 4 & 4 & 4 & 4\\
  &   &   & 3 &   &   &   &   & 3 & 3 & 3 & 3 & 3 & 3 & 3 & 3\\
  &   & 2 &   &   &   &   &   & 2 & 2 & 2 & 2 & 2 & 2 & 2 & 2\\
  & 1 &   &   &   &   &   &   & 1 & 1 & 1 & 1 & 1 & 1 & 1 & 1\\
0 &   &   &   &   &   &   &   & 0 & 0 & 0 & 0 & 0 & 0 & 0 & 0\\
\hline
\end{tabular}
\caption{Matrix representation of the staircase scheme
with $M=16$ and $B=8$ (entries indicate the buses connected to
memory modules).}
\end{center}
\end{table}

\vskip 1cm

\subsection*{Acknowledgments}  
The authors thank professors J.L.A. Yebra and E. Garriga for their valuable comments.
This research was supported by the {\em Ministerio de Econom\'{\i}a y Competitividad} and the {\em European Regional Development Fund} under project MTM2011-28800-C02-01, and by the Catalan Government
under project 2014SGR1147.


\begin{thebibliography}{99}

\bibitem{az98}
M. Aigner, G.M. Ziegler, {\it Proofs from THE BOOK}, Springer,
Berlin, 1998.


\bibitem{Ap77}
K. Apple, An attempt to understand the four color problem,
{\it J. Graph Theory} \textbf{1} (1977) 193--206.

\bibitem{ApHaKo77}
K. Apple, W. Haken, J. Koch,
Every planar map is four colorable,
{\it Illinois J. Math.} \textbf{21} (1977) 429--567.

\bibitem{b93}
N.L. Biggs, {\it Algebraic Graph Theory}, Cambridge University
Press, Cambridge, 1993.

\bibitem{blw76}
N.L. Biggs, E.K. Lloyd, R.J. Wilson, {\it Graph Theory: 1736-1936},
Claredon Press, Oxford, 1976.

\bibitem{Bo04}
B. Bollob\'{a}s, {\it Extremal Graph Theory},
Dover Publications, Mineola, N.Y., cop. 2004.

\bibitem{Bo90}
B. Bollob\'{a}s, {\it Graph Theory: An Introductory Course},
Springer, New York, 1979, 3rd corrected edition, 1990.

\bibitem{bcn89}
A.E. Brouwer, A.M. Cohen,  A. Neumaier, {\it Distance-Regular
Graphs}, Springer-Verlag, Berlin, 1989.

\bibitem{BrHa95}
A.E. Brouwer, W.H. Haemers, Association schemes,
in: R.L. Graham, et al. (eds.), \emph{Handbook of
Combinatorics}, Vol. 1--2, Elsevier, Amsterdam, 1995, 747--771.

\bibitem{Ca74}
L. Carroll, The Hunting of the Snark, Annotated by M. Gardner,
Penguin Books, New York, 1974.

\bibitem{Di97}
R. Diestel,
\emph{Graph Theory}.
Springer, New York, 1997.

\bibitem{ers66}
P. Erd\H{o}s, A. R\'{e}nyi, V. S\'{o}s, On a problem of graph
theory, {\it Studia Sci. Math.} \textbf{1} (1966) 215--235.

\bibitem{e1741}
L. Euler, Solutio problematis ad geometriam situs pertinentis, {\it
Commentarii academiae scientiarum petropolitanae} \textbf{8} (1741)
128--140.

\bibitem{e1752}
L. Euler, Elementa doctrine solidorum, {\it
Novi comm. acad. scientiarum imperialis petropolitanae} \textbf{4} (1752-1753)
109--160.

\bibitem{ff89}
J. F\`abrega, M.A. Fiol, Maximally connected digraphs, {\it J.
Graph Theory} {\bf 13} (1989) 657--668.

\bibitem{f91}
M.A. Fiol, A Boolean algebra approach to the construction of snarks,
in {\it Graph Theory, Combinatorics and Applications}, Y.~Alavi, G.~Chartrand, O.R.~Oellermann and
A.J.~Schwenk (eds.), Vol. 1, John~Wiley \& Sons, New York, 1991,
493--524.

\bibitem{f93}
M.A. Fiol, The connectivity of large digraphs and graphs, {\it J.
Graph Theory} {\bf 17} (1993) 31--45.

\bibitem{f95}
M.A. Fiol,
$c$-Critical graphs with maximum degree three, in:
{\it Graph Theory, Combinatorics, and Applications}
Y. Alavi and A.J. Schwenk, (eds.), Vol. 1, John~Wiley \& Sons, New York, 1995, 403--411.

\bibitem{f99}
M.A. Fiol, Algebraic characterizations of distance-regular graphs,
{\it Discrete Math.} \textbf{246} (2002), no.1-3, 111--129.

\bibitem{f&f84}
M.A. Fiol, M.L. Fiol,
Coloracions: un nou concepte dintre la teoria de coloraci\'{o} de grafs,
{\it L'Escaire} \textbf{11} (1984) 33--44 (in Catalan).

\bibitem{fg97}
M.A. Fiol, E. Garriga, From local adjacency polynomials to
locally pseudo-distance-regular graphs, {\it J. Combin. Theory Ser.
B} {\bf 71} (1997) 162--183.

\bibitem{fg99}
M.A. Fiol, E. Garriga, On the algebraic theory of
pseudo-distance-regularity around a set, {\it Linear Algebra Appl.}
{\bf 298} (1999), no. 1-3, 115--141.

\bibitem{fgy95}
M.A. Fiol, E. Garriga, J.L.A. Yebra, Avoiding monocoloured
triangles when colouring $K_n$, {\it Research Report}, Universitat
Polit\`{e}cnica de Catalunya, Barcelona, 1995.

\bibitem{fvyl84}
M.A. Fiol, M. Valero, J.L.A. Yebra, T. Lang, Reduced
interconnection networks based on the multiple bus for
multimicroprocessor systems, {\it International Journal of Mini and
Microcomputers} {\bf 6} (1984), no. 1, 4--9.

\bibitem{Ga76}
M. Gardner, Mathematical games: Snarks, Boojums and other conjectures related
to the four-color-map theorem, \emph{Sci. Amer.} \textbf{234} (1976) 126--130.

\bibitem{g93}
C.D. Godsil, {\it Algebraic Combinatorics}, Chapman and Hall, New
York, 1993.

\bibitem{GoRo01}
C.D. Godsil, G. Royle,
{\it Algebraic Graph Theory}, Springer-Verlag, New York, 2001.

\bibitem{h69}
F. Harary, {\it Graph Theory}, Addison-Wesley, Reading, MA, 1969.

\bibitem{h1890}
P.J. Heawood, Map Colour Theorems, \emph{Quart. J. Math.}
\textbf{24} (1890) 332--338.

\bibitem{Ho98}
P. Hoffman,
{\it The Man Who Loved Only Numbers: The Story of Paul Erd\H{o}s and the Search for Mathematical Truth}.
Hyperion, New York, 1998.

\bibitem{i75}
R. Isaacs,
Infinite families of nontrivial graphs which are not Tait colorable,
{\it Amer. Math. Monthly} {\bf 82} (1975), no. 3, 221--239.

\bibitem{i76}
R. Isaacs, Loupekhine's snarks: a bifamily of
non-Tait-colorable graphs, {\it Technical Report}, {\bf 263}, Dept.
of Math Sci., The Johns Hopkins University, Maryland, 1976.

\bibitem{Ku30}
C. Kuratowski,
Sur le probl\`{e}me des courbes gauches en topologie,
{\it Fund. Math.} {\bf 15} (1930) 217--283.

\bibitem{lvf83}
T. Lang, M. Valero, M.A. Fiol, Reduction of connections for
multibus organization, {\it IEEE Trans. Comput.} {\bf C-32} (1983),
no. 8, 707--716.

\bibitem{Pe1898}
J. Petersen, Sur le th\'{e}or\`{e}me de Tait,
\emph{L'Interm\'{e}diaire des Math\'{e}maticiens} \textbf{5} (1898) 225--227.

\bibitem{RaTo70}
H. Rademacher, O. Toeplitz, {\it The Enjoyment of Math},
Princenton University Press, Princenton, New Jersey, 1957.

\bibitem{r06}
S.P. Radziszowski, Small Ramsey numbers, {\it Electron. J. Combin.}
{\bf 1} (2006) Dynamical Survey 1.

\bibitem{RoSaSeTh97}
N. Robertson, D.P. Sanders, P.D. Seymour, R. Thomas,  The four colour theorem, {\it
J. Combin. Theory Ser. B.} {\bf 70} (1997) 2--44.

\bibitem{Ta1880}
P.G. Tait, Remarks on the colouring of maps, {\it
Proc. Roy. Soc. Edimburgh} {\bf 10} (1880) 501--503, 729.

\bibitem{Tho90}
C. Thomassen,
A link between the Jordan curve theorem and the Kuratowski planarity criterion, {\it Amer. Math. Monthly} {\bf 97} (1990) 216--218.

\bibitem{vl99}
J.H. van Lint, {\it Introduction to Coding Theory}, third edition,
Springer, Berlin, 1999.

\bibitem{w00}
D.B. West, {\it Introduction to Graph Theory}, Prentice-Hall, Englewood Cliffs, NJ, 2nd edition, 2000.
\end{thebibliography}
\end{document}